\documentstyle[12pt]{article}

\begin{document}

\newcommand{\supno}{\stackrel{\supset}{\not=}}
\newcommand{\subno}{\stackrel{\subset}{\not=}}
\newcommand{\zbz}{Z_{B_0}}
\newcommand{\zbe}{Z_{B_1}}
\newcommand{\zbk}{Z_{B_2}}
\newcommand{\zm}{Z_{B_m}}
\newcommand{\zme}{Z_{B_{m-1}}}
\newcommand{\zmk}{Z_{B_{m-2}}}
\newcommand{\zi}{Z_{B_i}}
\newcommand{\zj}{Z_{B_j}}
\newcommand{\bc}{{\bf C}}
\newcommand{\bz}{{\bf Z}}
\newcommand{\g}{A}				
\newcommand{\zn}{Z_{num}}				
\newcommand{\zk}{Z_K}	
\newcommand{\co}{{\cal O}}
\newcommand{\x}{$(X,p)$\ }		
\newcommand{\cl}{{\cal L}}
\newcommand{\ci}{{\cal I}_{D'}}
\newcommand{\y}{S. S.-T. Yau\ }

\noindent 
{\large\bf ``Weakly''  Elliptic Gorenstein Singularities of Surfaces}\\

\noindent 
{\bf  Andr\'as N\'emethi}\footnote{Partially supported by 
NSF Grant No. DMS-9622724}\\

\noindent 
Department of Mathematics, The Ohio State University, Columbus OH 43210, USA,\\
(e-mail: nemethi@math.ohio-state.edu)\\

\noindent {\bf 1.\ Introduction.}\\

Let $p$ be a singularity of a normal two--dimensional analytic surface $X$.
In general
$\pi:M\to X$ will denote a resolution of $(X,p)$, and $\g=\pi^{-1}(p)$. 
It is well--known that 
the dual resolution graph of the germ \x determines the 
topology of the germ 
$(X,p)$. Indeed, the (real) link  of \x can be reconstructed from 
the resolution graph as a plumbing. Conversely, by a result of W. Neumann
\cite{Neumann}, the topology of \x determines  the (minimal) resolution
graph. Therefore, if an invariant of \x depends only on the 
dual resolution  graph, we say that it is a ``topological invariant''. 
For example, the self--intersection $\zn^2$ of Artin's (fundamental) 
cycle $\zn$, or the Euler--characteristic
$\chi(D)$ of any cycle $D$ supported by $\g$, can be determined from 
the graph. 

Now, it is fascinating to investigate  if an invariant, 
a priori defined from the analytic structure of $(X,p)$, is topological
or not. In this article, we ask this question for
the geometric genus $p_g=h^1(M,\co_M)$ and the Hilbert--Samuel  function
of $(X,p)$, in particular
for the 
multiplicity $mult$\x and the embedding 
dimension $emb\,dim$\x$=\dim m_p/m^2_p$.
In general, these invariants are not topological, but if we restrict our
study to some special classes, then they can be determined from the graph.

The first result of this type was obtained by M. Artin \cite{Artin62,Artin66}
for rational singularities.
He proved that they can be characterized topologically:
$$\mbox{\x is rational, i.e. $p_g=0$} \Leftrightarrow
\chi(\zn)=1\Leftrightarrow \min_{D>0}\chi(D)\geq 1.$$
Moreover, for these singularities, 
the Hilbert--Samuel function (hence $mult$\x and $emb\,dim$\x too) can also be 
computed from the resolution graph. 

The elliptic singularities were introduced by P. Wagreich in \cite{Wagreich}
(cf. ``Terminology'' at the end of the introduction). 
They are defined topologically:
$$\mbox{\x is elliptic}\Leftrightarrow\chi(\zn)=0\Leftrightarrow \min_{D>0}
\chi(D)=0.$$

The class of elliptic singularities
contains all the singularities with $p_g=1$, and all the 
Gorenstein singularities with $p_g=2$; but
an  elliptic singularity can have arbitrary high geometric genus
(see e.g. the examples after Theorem C).
The next step in the above program was obtained by H. Laufer \cite{Laufer77}.
He proved that the Gorenstein singularities with $p_g=1$ (he called them 
minimally elliptic singularities) can be characterized topologically
(cf. 2.7), and for these singularities
all the above (a priori) analytical invariants are topological.

Moreover, he noticed that singularities with $p_g=1$ (without Gorenstein
assumption), or even Gorenstein singularities with $p_g=2$ (or $p_g\geq 2$)
cannot be characterized topologically. In this second case, one can easily
construct pairs of hypersurface singularities with the same resolution graph
but different $p_g$ (cf. 2.22). But in all these examples (known by the
author) either there is a non-rational exceptional divisor in the resolution,
or the graph is not a tree, i.e. $H^1(\g,\bz)\not=0$. (Notice that
$H^1(\g,\bz)=0$ if and only if the link of \x is a rational homology sphere.)

\vspace{1mm}

{\em The main message of the present paper is that for Gorenstein singularities
with $H^1(\g,\bz)=0$ the Artin--Laufer program can be continued.
Here we give the complete answer in the case of elliptic singularities.}

\vspace{1mm}

Elliptic singularities were intensively studied by Wagreich \cite{Wagreich},
 Laufer \cite{Laufer77}, \y 
\cite{Yau4,Yau6,Yau7,Yau2,Yau3,Yau5,Yau1} and others. 
In Yau's papers, the fundamental topological invariant is the 
 ``elliptic sequence''.
In the numerical Gorenstein case its definition is the following.
Let $\pi:M\to X$ be the {\em minimal} resolution of $(X,p)$, and 
$\g=\pi^{-1}(p)$ as
above. Let $\zk$ be the canonical cycle (cf. \S 2). 
The elliptic sequence consists
of a  sequence $\{\zj\}_{j=0}^m$, where $\zj$ is Artin's
(fundamental) 
cycle of $B_j\subset \g$. We define $\{B_j\}_j$ inductively as follows.
For $j=0$ take $B_0=\g$ hence $\zbz=\zn$. Then  $\zk\geq \zbz$. If $\zk>\zbz$
then we set
$B_1:=|\zk-\zbz|$. Similarly, 
if $B_i$ is already  defined for any $i\leq j$,
then $\zk\geq \zbz+\cdots+\zj$
(for details, see 2.10). If the inequality is 
strict then we define $B_{j+1}:=|\zk-\zbz-\cdots-\zj|$, otherwise we stop.
The length of the elliptic sequence $\{\zj\}_{j=0}^m$ is
$m+1$. The case $m=0$ corresponds exactly to the minimally elliptic
singularities of Laufer \cite{Laufer77}. 

\y proved that for a numerical Gorenstein elliptic
singularity $p_g\leq m+1$ (\cite[(3.9)]{Yau1}, cf. also 2.19).
Particular examples show that strict inequality
can occur. 

In the first sections we give several characterizations of the ``extremal
property'' $p_g=m+1$. The most important characterization,
from the point of view of the present manuscript, is the following:

\vspace{2mm}

\noindent {\bf Theorem A.}\ 
{\em Let  \x be  an elliptic numerical Gorenstein singularity and set
$C_j'=\sum_{i\geq j}Z_{B_i}$. Then
$p_g=m+1$ if and only if 
the ``obstruction line bundles''
 $\co_{C'_j}(-Z_{B_{j-1}})$ are trivial
for all $1\leq j\leq m$.}

\vspace{2mm}

\noindent Now, the main technical  result of the paper (see section 4) says:

\vspace{2mm}

\noindent {\bf Theorem B.}\ {\em 
Let \x be  an elliptic Gorenstein  singularity. 
Fix an integer $1\leq k\leq m$ and 
assume that for $k+1\leq j\leq m$ the ``obstruction line bundles''
$\co_{C'_j}(-Z_{B_{j-1}})$ are trivial. (If $k=m$
then this assumption is vacuous.) Then $\co_{C'_k}(-Z_{B_{k-1}})$ 
has finite order. (Actually, its order is not greater  than $k!$.)}

\vspace{2mm}

\noindent Now, 
 if $H^1(\g,\bz)=0$ then $Pic(D)$ is torsion free for any positive cycle
$D$. Hence, we obtain the main result of the paper:

\vspace{2mm}

\noindent 
{\bf Theorem C.}\ {\em Assume that \x is an elliptic Gorenstein
singularity with  $H^1(\g,\bz)=0$. 
Then  $p_g$ is a topological invariant. In fact $p_g=m+1=$ the length of 
the elliptic sequence in  the minimal resolution of $(X,p)$. }

\vspace{2mm}

\noindent 
Actually, one of the starting points of our investigation was S. S.-T. Yau's
result in \cite{Yau3}, which says that 
{\em hypersurface} singularities with $H^1(\g,\bz)=0$ and $p_g=2$
 satisfies $m=1$ (cf. Remark 4.14). 
(In Yau's terminology, a numerical Gorenstein elliptic singularity 
with $p_g=m+1$ is called ``maximally elliptic''.)

\vspace{2mm}

\noindent {\bf Examples.}\ \cite{Wagreich,Yau1}\ Set $(X_i,0)\subset (
{\bf C}^3,0)$ given by 
$(X_1,0)=\{z^2+y^3+x^{9+6m}=0\}$,
$(X_2,0)=\{z^2+y^3+x^{11+6m}=0\}$,
$(X_3,0)=\{z^3+y^3+x^{4+3m}=0\}$,
$(X_4,0)=\{z^3+y^3+x^{5+3m}=0\}$ (where $m\geq 0$).
Then, in all these cases, $H^1(\g_i,\bz)=0$, $(X_i,0)$ is elliptic,
and $p_g=m+1=$ the length of the elliptic sequence ($1\leq i\leq 4$)
(cf. also with 6.4.c).

\vspace{2mm}

\noindent 
In sections 5 and 6,
 we generalize Laufer's results about the minimally elliptic
singularities \cite{Laufer77} (which corresponds to $m=0$). Namely, we prove:

\vspace{2mm}

\noindent {\bf Theorem D.}\  {\em Assume that \x is an elliptic
Gorenstein singularity with $p_g=m+1$.
(If $m=0$ or $H^1(\g,\bz)=0$ then the last assumption is satisfied.)
Let $\zn$ be Artin's cycle in the minimal resolution $(M,\g)\to (X,p)$. Then:

(a)\ If $\zn^2\leq -2$, then $m_p\co_M=\co(-\zn)$, hence $mult(X,p)=-\zn^2$;

(b)\ 
If $\zn^2=-1$, then $m_p\co_M=m_Q\co(-\zn)$ for some smooth point $Q$
of $\g$, and $mult(X,p)=2$;

(c)\ $\mbox{emb dim\x}=\max(3,-\zn^2)$;

(d)\
$H^0(M,\co(-k\zn))=m_p^k$ \ for any $k\geq 0$ provided that $\zn^2\leq -3$;

(e)\ 
$\dim\co_{(X,p)}/m_p^k=\chi(\co_{k\zn})+1$, and 
$\dim m_p^k/m_p^{k+1}=-k\zn^2$ \  for any $k\geq 1$ provided that 
$\zn^2\leq -3$;
(for the cases $\zn^2=-1$ or $-2$, see (6.4)). }

\vspace{2mm}

Some partial results in this direction were obtained by \y 
in his series of papers (we will cite them in the body of the paper
at the corresponding places). 
In order to have a self--contained presentation,
we reprove those facts what we use from his work. 
In some cases we follow Yau's original arguments, 
in other  cases we give  different  proofs (original ones or arguments in 
the spirit of \cite{MR}). The article of M. Reid \cite{MR} was 
extremely helpful for the author. Actually, the proofs 
in sections 5 and 6 have their origins in \cite{MR}.
One of the ideas of the proof of Theorem B was borrowed from J. Wahl's paper
\cite{Whal}.
On the other hand, we emphasize that almost all the classical arguments,
used in the case of rational or  minimally elliptic singularities, and  based
on some vanishing theorems or on the numerically 
1-- or 2--connectivity of $\zn$, in our 
general situation fail, and we had to replace them by different arguments.

Some of the results of the present article
(especially, the multiplicity and embedding dimension
computations)  can be compared with the results
of U. Karras and J. Stevens proved for Kodaira, respectively Kulikov
singularities. For details, see \cite{K1,K2} and \cite{St1,St2}.

Finally, we notice that Theorem D is not true if $\chi(\zn)<0$.
E.g. for $(X,0)=\{x^3+y^4+z^7=0\}$ one can verify that $\chi(\zn)=-1$, 
$H^1(\g,\bz)=0$, and  $\zn^2=-2$, but $mult(X,0)=3$.
\\

\noindent \ {\bf Terminology:}\
Singularities characterized by $\chi(\zn)=0$ sometimes are called 
``weakly'' elliptic. By this terminology one wants to emphasize the 
difference between these singularities and the ``strongly'' elliptic
singularities, defined by $p_g=1$; cf. also Yau's papers.
Notice also that in the terminology
of M. Reid \cite{MR}, ``elliptic'' means ``minimally elliptic'' in the sense of
Laufer (cf. 2.7). In this article  we will adopt the terminology
used by Wagreich and Laufer.\\

The author wishes to thank \'Ecole Polytechnique at Palaiseau,
University of Nice and University of Nantes
(especially Professors C. Sabbah, M. Merle and F. Elzein) for their hospitality
and excellent working atmosphere.\\

\noindent {\bf 2.\ Preliminaries.}\\

We fix a normal surface singularity $(X,p)$. Let $\pi:M\to X$ be its
{\em minimal} resolution, and let $\g=\cup_i\g_i$ be the decomposition of 
the exceptional set $\g$ into irreducible components.
Let $K$ be the canonical divisor on $M$, hence:
$$(2.1)\hspace{3cm}
\g_iK=-\g_i^2-2+2g_i+2\delta_i\hspace{1cm}\mbox{for all $i$,}\hspace{2cm}$$
where $g_i$ is the genus of $\g_i$, and $\delta_i$ is the sum of all 
delta-invariants of the singular points (the ``number of nodes and cusps'')
on $\g_i$ \cite{Serre}. 

In this paper all the cycles will be integer combinations of the $\g_i$'s. If
$D=\sum_in_i\g_i$, we write $n_i=m_{\g_i}(D)$. $|D|$ denotes the 
support of $D$.

We denote Artin's numerical (fundamental) cycle by $\zn$, i.e. $\zn$ is the 
minimal positive cycle 
$Z$ with $Z\g_i\leq 0$ for all $i$ \cite{Artin62,Artin66}. 
\x is called numerical Gorenstein if there is a cycle $Z_K$ 
(called canonical cycle) which satisfies 
$\zk\g_i=-K\g_i$ for all $i$. If $\zk$ exists, then  
(by 2.1) $Z_K$ is trivial if and only if \x is a Du Val singularity. 
Otherwise $m_{\g_i}(Z_K)>0$ for all $\g_i$ (cf. \cite{Laufer87}, page 490).
Moreover, if $\zk\not=0$, then from (2.1) 
one has $Z_K\g_i\leq 0$ for all $i$, hence $\zn\leq\zk$.

By a (Kodaira type) vanishing theorem (cf. \cite{MR}, Exercise 15, page 119),
for any  cycle $D\geq 0$ which
satisfies $\g_iD\leq 0$ for all $\g_i$, and any line bundle $\cl$ with
$deg_{\g_i}\cl\geq K\g_i$ for all $i$:
$$(2.2)\hspace{.5cm}
\left\{\begin{array}{l}
h^1(M,\co(-D)\otimes \cl)=0,
\mbox{hence in the numerical Gorenstein case:} \\
H^0(M,\co(-D))\to H^0(\co_{\zk}(-D))\ \mbox{is onto.}
\end{array}\right.$$

For any cycle $D$, we denote the Euler--characteristic $h^0(\co_D)-h^1
(\co_D)$ by $\chi(D)$. By Riemann--Roch theorem 
$\chi(D)=-D(D+K)/2$. By \cite{Laufer77}, (2.6):
$$(2.3)\hspace{2cm} 
H^0(\co_{\zn})=\bc=\{\mbox{constants}\},  \ 
\mbox{hence}\ \ \chi(\zn)\leq 1.\hspace{2cm} $$

\noindent By a result of Artin \cite{Artin62,Artin66}, the following facts 
are equivalent:
$$(2.4)\hspace{1cm}
 \mbox{\x is rational (i.e. $p_g=0$)}\Leftrightarrow
\chi(\zn)=1\Leftrightarrow \min_{D>0}\chi(D)\geq 1.$$
The dual resolution graph of a rational singularity is a tree, and all the 
vertices corresponds to smooth rational curves.

In \cite{Wagreich}, Wagreich introduced the elliptic singularities. They are
defined by the property $\min_{D>0}\chi(D)=0$. 
In particular, by (2.3) and (2.4), $\chi(\zn)=0$. 
The inverse implication  is also true, see e.g. \cite{Laufer77} (4.2), hence:
$$(2.5)\hspace{2cm}
\mbox{\x is elliptic}\Leftrightarrow\chi(\zn)=0\Leftrightarrow \min_{D>0}
\chi(D)=0.\hspace{2cm}$$
Elliptic singularities include all the singularities with $p_g=1$ and all
the Gorenstein singularities with $p_g=2$ (cf. the proof of  4.13);
but elliptic singularities can have arbitrary high geometric genus.

Following  Laufer (Definition 3.1 \cite{Laufer77}), we say  that a cycle $E>0$
is minimally elliptic if $\chi(E)=0$ and $\chi(D)>0$ for all cycles 
$0<D<E$. Laufer in \cite{Laufer77} (3.2) proved that if $\chi(\zn)=0$
then there exists a unique minimally elliptic cycle $E$. Hence,
in the elliptic numerical Gorenstein case:
$$(2.6)\hspace{5cm} E\leq \zn\leq \zk.\hspace{5cm}$$
The minimally elliptic singularities are characterized by the following
equivalent properties 
(cf. Laufer's paper \cite{Laufer77})
(notice that (b), (c) and (d) are topological properties; (d) provides the 
name of the singularity):
$$(2.7)\hspace{1cm} \left\{\begin{array}{l}
(a)\ \mbox{\x is Gorenstein and $p_g=1$};\\
(b)\ \mbox{$\zn$ is a minimally elliptic cycle};\\
(c)\ \mbox{\x is numerical Gorenstein with $\zk=\zn$};\\
(d)\ \mbox{$\chi(\zn)=0$ and any connected proper subvariety of $\g$}\\
\ \ \ \ \mbox{supports a rational singularity.}\end{array}\right.$$

\noindent 2.8.\ If \x is elliptic, then the topology of the irreducible
exceptional divisors $\g_i$ cannot be very complicated. Wagreich in
\cite{Wagreich}, page 428, proved that precisely one of the followings hold:
(a) Precisely one component $\g_{i_0}$ satisfies $\chi(\g_{i_0})=0$
(which is either a smooth elliptic curve or a rational curve with 
$\delta_{i_0}=1$)
and the other irreducible exceptional divisors are smooth rational curves.
(In this case $E=\g_{i_0}$). Or (b): all the exceptional divisors are smooth rational curves.

Actually, from the uniqueness of the minimally elliptic cycle it follows that
all the connected components of $\g\setminus |E|$ support the exceptional
set of rational singularities, hence the dual graph of \x can be obtained 
from the dual graph of $|E|$ by gluing trees whose vertices correspond to
smooth rational curves. In particular, the restriction map
$$(2.9)\hspace{2cm} 
H^1(\g,\bz)\to H^1(|E|,\bz)\ \mbox{is an isomorphism.}\hspace{2cm} $$
In a series of papers \y investigated the properties of 
elliptic singularities (he called them ``weakly elliptic''), cf.
\cite{Yau4,Yau6,Yau7,Yau2,Yau3,Yau5,Yau1}. His main tools are the 
``computation sequence'' 
introduced by Laufer in \cite{Laufer72} and \cite{Laufer77},
and the  ``elliptic sequence''. 
The interested reader can find in Yau's  papers the definition of the 
elliptic sequence  for an arbitrary elliptic singularity. In the 
non--numerical Gorenstein case the sequence has a lot of anomalies (see e.g.
some  examples  on the page  881 of \cite{Yau5}), but in the numerical
Gorenstein case it is a powerful tool. Since the definition in this later case
is easier and more natural, and it is sufficient for our goals, we adopt the general definition for this situation.\\
2.10.\ {\bf The construction of the elliptic sequence.}\ Consider the minimal resolution of an elliptic numerical Gorenstein singularity. The elliptic sequence consists
of the sequence $\{\zj\}_{j=0}^m$, where $\zj$ is the numerical (fundamental) 
cycle of $B_j\subset \g$. We define $\{B_j\}_j$ inductively as follows.
For $j=0$ take $B_0=\g$ hence $\zbz=\zn$. By (2.6) $E\leq \zn\leq \zk$. 

If $\zbz<\zk$, we define $B_1:=|\zk-\zn|$. 
Then $B_1$ is connected. Indeed, assume that $B_1$ has more connected 
components $\{B_{1,t}\}_t$; then write $\zk-\zn=\sum_t D_t$ with 
$|D_t|=B_{1,t}$. Then the vanishing  $\chi(\zk-\zn)=0$ and (2.5) imply that 
$\chi(D_t)=0$ for all $t$. Then by \cite{Laufer77} (3.2),
each $B_{1,t}$ supports  a minimally elliptic cycle. This 
contradicts  the uniqueness of the minimally elliptic cycle  $E$.
Now, using again 
$\chi(\zn)=0$, one gets $\zn(\zk-\zn)=0$. Since $\zn\g_i\leq 0$ for all 
$\g_i$, it follows that for any $\g_i\subset B_1$ the equality $\zn\g_i=0$
holds. In particular, by the non-degeneracy of the intersection form,
$B_1\not= B_0$. Moreover, for any $\g_i\subset B_1$ one has:
$\g_i(\zk-\zbz)=\g_i\zk$, 
hence $B_1$ supports a numerical Gorenstein singularity
with canonical cycle $\zk-\zbz$. This singularity is again elliptic. Indeed,
$\min_{D\subset B_1}\chi(D)\geq \min_{D\subset B_0}\chi(D)=0$, but
$\chi(\zk-\zn)=0$, hence $\min_{D\subset B_1}\chi(D)=0$. In particular,
(cf. 2.6), $E\leq \zbe\leq \zk-\zbz$. 

Now we repeat the above arguments. 
If $\zbe<\zk-\zbz$, we define $B_2:=|\zk-\zbz-\zbe|$ and we verify 
that it supports an elliptic singularity with canonical cycle $\zk-\zbz-\zbe$ 
and for any $\g_i\subset B_2$ the vanishing $\g_i\zbe=0$ holds. If we repeat this procedure, after a finite step we will obtain $\zm=\zk-\zbz-\cdots-\zme$,
i.e. the numerical cycle and the canonical cycle of $B_m$ 
coincides, hence by (2.7) $B_m$ supports a minimally elliptic singularity
with $E=\zm$. In particular $B_m=|E|$. 

If we contract the connected exceptional curve  $B_j\subset M$, we obtain a
unique singular point; this will be denoted by $(M/B_j,p_j)$ 
($0\leq j\leq m$). It is convenient to introduce the notations
$C_t=\sum_{i=0}^t\zi$ and $C'_t=\sum_{i=t}^m\zi$ $(0\leq t\leq m)$ too.

\vspace{1mm}

\noindent 
2.11.\ {\bf Definition/first properties of the elliptic sequence:}\

{\em 
(a)\ $B_0=\g,\ B_1=|\zk-\zbz|,\ B_2=|\zk-\zbz-\zbe|,\ \ldots, \ B_m=|E|; $
each $B_j$ is connected and the inclusions $B_{j+1}\subset B_j$ are strict.
Moreover, $\zn=\zbz\supset\zbe\supset \cdots\supset \zm=E$.

(b)\ If $\g_i\subset B_{j+1}$ then $\g_i\zj=0$ for all $i$ and $j$. 
In particular, $\zi\zj=0$ for all $0\leq i<j\leq m$.

(c)\ $Z_K=\sum_{i=0}^m\zi$ (cf. also \cite[(3.7)]{Yau1}).

(d)\ $\g_iC'_t=\g_i\zk$ for any $\g_i\subset |C'_t|$. In other words,
$C'_t$ is the canonical cycle of $|C'_t|=B_t$ (i.e. of $(M/B_t,p_t)$).

(e)\ For any $\g_i\subset \g$ the inequality $\g_i\cdot C_t\leq 0$ holds.}

\vspace{1mm}

\noindent 
{\em Proof.}\ $a,b,c,d$ follows from the above construction. The proof of
$(e)$ is as follows.
If $\g_i\subset B_t$ then $\g_i\zj\leq 0$ for any $j\leq t$
by the definition of the numerical cycle, hence $\g_iC_t\leq 0$.
If $\g_i\not\subset B_t$ then $\g_iC_t=\g_i(\zk-C'_{t+1})$. Now, $\g_i\zk\leq
0$ (by the minimality of the resolution) and $\g_iC'_{t+1}\geq 0$
(because $|C'_{t+1}|\subset B_t$). \ \ \ $\Box$

\vspace{2mm}

\noindent 2.12.\ The next result shows that the properties $(d)$ and $(e)$
characterize the elliptic sequence.\\
2.13. \ {\bf Lemma.}\ 
{\em  Assume that \x is a  numerical Gorenstein elliptic singularity.

(a)\  If a
cycle $Z\geq 0$ satisfies $\g_i(Z-\zk)\geq 0$ for all $\g_i\subset |Z|$, then 
in fact $\g_i(Z-\zk)=0$ for all $\g_i\subset |Z|$, and 
$Z\in \{C'_0,C'_1,\ldots,C'_m,0\}$.

(b)\ If a cycle $0\leq Z\leq \zk$ satisfies  $\g_iZ\leq 0$ for all 
$\g_i\subset \g$, then $Z\in\{0,C_0,C_1,\ldots,C_m\}$. }

\vspace{1mm}

\noindent {\em Proof.}\ 
If $\g_i(Z-\zk)\geq 0$, then $Z(Z-\zk)\geq 0$, hence $\chi(Z)\leq 0$. This 
together with (2.5) guarantees that $\chi(Z)=0$ and $\g_i(Z-\zk)=0$
for all $\g_i\subset |Z|$. Now, 
if $|Z|=B_0$, then by the non--degeneracy of the intersection form $Z=\zk$.
If $|Z|\not= B_0$, then for $\g_i\subset |Z|$ obviously $\g_i(\zk-Z)=0$;
for $\g_i\not\subset |Z|$ one has $\g_i(\zk-Z)\leq 0$ (because $\g_i\zk\leq 
0$ and $\g_iZ\geq 0$). Hence  $\g_i(\zk-Z)\leq 0$ for any $\g_i$, therefore
$\zk-Z\geq \zbz$, hence $Z\leq \zk-\zbz$. In particular, $|Z|\subset B_1$. 
If $|Z|=B_1$, then $Z$ is  the canonical cycle $C'_1$ of $B_1$, 
otherwise, by the same argument as above $Z\leq \zk-\zbz-\zbe$, hence
$|Z|\subset B_2$. Continuing the precess, (a) follows.
For the second part, apply (a) for $\zk-Z$. \ \ \ $\Box$\\

\noindent Since $E\leq \zj\leq \zn$,
$$(2.14) \hspace{1cm}
h^1(\co_{\zj})=1\ \ \mbox{hence}\ \ 
\chi(\zj)=\chi(C_j)=\chi(C'_j)=0
\ \ \mbox{for all}\ \ 0\leq j\leq m.$$

\noindent 2.15.\  Since any $\zi$ is a numerical fundamental cycle, 
by \cite{Laufer72} (proof of 4.1) or \cite{Laufer77} (cf. also \cite{Yau1}),
there is a ``computation sequence'' which can 
start with  any of the irreducible exceptional divisors $\g_j
\subset B_i$ and ends with
$\zi$. This means that there exists a sequence of cycles 
$Z_0,\ldots , Z_k$ with $Z_0=\g_j, Z_k=\zi,\
Z_{l+1}=Z_l+\g_{i_l}$ and $\g_{i_l}Z_l>0$ for $0\leq l\leq k-1$. Moreover, 
for any  two integers $0\leq i<j\leq m$, since $\zi$ is a 
numerical cycle  with $\zj<\zi$,  there is a computation 
sequence which starts with $\zj$ and ends with $\zi$. More precisely, there is 
a sequence  $Z_0,\ldots , Z_k$ with $Z_0=\zj, Z_k=\zi,\
Z_{l+1}=Z_l+\g_{i_l}$ and $\g_{i_l}Z_l>0$ for $0\leq l\leq k-1$. 
Actually in this second case, since $h^1(\co_{\zi})=
h^1(\co_{\zj})=1$ (cf. 2.14),  for any  $0\leq l\leq k-1$
the curve $\g_{i_l}$ is smooth rational and $\g_{i_l}Z_l=1$
(see also \cite{Laufer77} (2.7)). We will call a 
sequence like this a ``computation sequence which connects $\zj$ and $\zi$''.

Let $D$ be a cycle with $D\g_k=0$ for any $\g_k\subset B_i$, and consider 
a computation sequence $\{Z_l\}_l$ which connects $\zj$ and $\zi$
(where $j>i$). In the 
exact sequence
$0\to \co_{\g_{i_l}}(-D-Z_l)\to \co_{D+Z_{l+1}}\to \co_{D+Z_l}\to 0$
$(0\leq l\leq k-1)$ the Chern number is 
$\g_{i_l}(-D-Z_l)=-1$, hence the natural map  
$$(2.16)\hspace{1cm}
H^k(\co_{D+\zi})\to H^k(\co_{D+\zj}) \  (k=1,2)\ \mbox{is an isomorphism}. $$
For example,
if $D=C_{i-1}$ and $j=m$, then 
$$(2.17)\hspace{5mm}
H^k(\co_{C_i})\to H^k(\co_{C_{i-1}+E})\ (1\leq i\leq m;\  k=1,2)\ 
 \ \mbox{is an isomorphism}.$$
Now consider  the exact sequence 
$$(2.18)\hspace{2cm}
0\to \co_E(-C_{i-1})\to \co_{C_{i-1}+E}\to \co_{C_{i-1}}\to 0.\hspace{2.5cm}$$
Since $E$ is 2--connected (see e.g. \cite[(4.21)]{MR}), 
$\co_E(-C_{i-1})$ is trivial if and only if
it has a non-zero section (see e.g. \cite{MR}, page 82), 
and in this case $h^k(\co_E(-C_{i-1}))=1$ for $k=1,2$.
Therefore (2.17) and (2.18) imply 
$h^1(\co_{C_{i-1}})\leq 
h^1(\co_{C_{i}})\leq h^1(\co_{C_{i-1}})+1$ for $1\leq i\leq m$. Since 
$h^1(\co_{C_0})=1$, one obtains \y's result
(\cite[(3.9)]{Yau1}):
$$(2.19)\ \ h^1(\co_{C_i})\leq i+1\ (\mbox{for}\ 0\leq i\leq m),
 \ \mbox{hence}\ p_g=h^1(\co_{C_m})\leq m+1.$$

The next theorem gives some 
characterizations  of  the ``extremal property''  $p_g=m+1$. 

\vspace{2mm}

\noindent 
2.20.\ {\bf First characterization of $p_g=m+1$.}\ 
{\em Assume that \x is a  numerical Gorenstein elliptic singularity. Then 
the following facts are equivalent:

(a)\ $p_g=m+1$;

(b)\ $h^1(\co_{C_i})= i+1$\ \  for all\ $0\leq i\leq m$;

(c)\ The line bundles $\co_E(-C_j)$ are trivial (in $Pic(E)$), and 
$H^0(\co_{C_j+E})\to H^0(\co_{C_j})$ is surjective for any $0\leq j\leq m-1
$;

(d)\ $h^1(\co(-C_j))=m-j$ for all $0\leq j\leq m$;

(e)\ $ h^1(\co_{C'_j})=m-j+1$ \ for all $0\leq j\leq m$;

(f)\ For any $0\leq j\leq m-1$, there exists $f_j\in H^0(M,\co(-C_j))$, such 
that for any $\g_l\subset  B_{j+1}$ the vanishing order of $f_j$ on $\g_l$ is
exactly $m_{\g_l}(C_j)$;

(g)\ The line bundles $\co_{C'_{j+1}}(-C_j)$ are trivial for $0\leq j\leq m-1$;

(h)\ The line bundles $\co_{C'_{j+1}}(-\zj)$ are trivial for $0\leq j\leq m-1$.
}\\
{\em Proof.}\ $a\Leftrightarrow b \Leftrightarrow c$ is just a reformulation 
of the above discussion (2.17-18-19).  (Actually $a\Leftrightarrow c$
can be found in  \cite{Yau7} (2.4).) 

$d\Rightarrow a$ is easy: for $j=0$ 
one has $m=h^1(\co(-\zn))\stackrel{(2.3)}{=}
h^1(\co)-h^1(\co_{\zn})\stackrel{(2.5)}{=}p_g-1$. 
Conversely, for  $(a,b,c)\Rightarrow d$, notice that 
$H^0(\co)\to H^0(\co_{C_m})$ is onto by (2.2), and $H^0(\co_{C_m})\to 
H^0(\co_{C_j})$ is onto  by part (c) and (2.17) applied several times.
Therefore $H^0(\co)\to H^0(\co_{C_j})$ is onto, hence
$h^1(\co(-C_j))=h^1(\co)-h^1(\co_{C_j})=m+1-j-1$.

$e\Rightarrow a$ again is trivial (take $j=0$). In order to prove 
 $a\Rightarrow e$, 
we construct a sequence $\{Z_l\}_l$. We start with $Z_0=\zk$.
If $Z_l>0$ is already constructed, then we define $Z_{l+1}=Z_l-\g_{i_l}$ as 
follows. If 
there is at least one $\g_i\subset |Z_l|$ with $\g_i(Z_l-\zk)<0$
then take for $\g_{i_l}$ one of these $\g_i$'s. Otherwise $\g_{i_l}$ is
an arbitrary $\g_i\subset |Z_l|$. In this second case, by (2.13),
$\g_{i_l}(Z_l-\zk)=0$ and $Z_l\in\{C'_0,\ldots,C'_m\}$. If $Z_l=0$ then we stop.

Now, from the exact sequence $0\to \co_{\g_{i_l}}(-Z_{l+1})\to \co_{Z_l}
\to \co_{Z_{l+1}}\to 0$ and Serre duality
one gets that $h^1(\co_{Z_l})=h^1(\co_{Z_{l+1}})$
in the first case, and in second case
$h^1(\co_{Z_{l+1}})=h^1(\co_{Z_{l}})-\epsilon_l$, where $\epsilon_l\in\{0,1\}$.
Since $h^1$ must drop exactly $m+1$ times, we obtain that in the above sequence
we reach all the cycles $C'_j$ and at every time $\epsilon_l=1$. 

$a\Rightarrow f$ is proved in  \cite{Yau1} (3.13). 
We present  a short proof of it.
Fix an arbitrary $\g_k\subset B_{j+1}$. 
Consider a computation sequence $\{Z_l\}_l$ of $Z_{B_{j+1}}$ which starts
with $\g_k$ (cf. 2.15). Since the Chern numbers $\g_{i_l}(-Z_l-C_j)<0$,
using the sequences 
$0\to \co_{\g_{i_l}}(-Z_{l}-C_j)\to 
\co_{Z_{l+1}}(-C_j)\to \co_{Z_{l}}(-C_j)\to  0$, we obtain that 
$\beta: H^0(\co_{Z_{B_{j+1}}}(-C_j))\to H^0(\co_{\g_k}(-C_j))$ is injective.
Now, $H^0(\co_{\g_k}(-C_j))=\bc$. Indeed, 
if $\g_k\approx {\bf P}^1$ then 
$H^0(\co_{\g_k}(-C_j))=H^0(\co_{{\bf P}^1})$. If $\g_k\not\approx 
{\bf P}^1$ then $\g_k=E$ (by 2.8), and $\co_{\g_k}(-C_j)$ is trivial by ($c$).
Now,  considers the diagram:

\begin{picture}(400,100)(0,-10)
\put(0,75){\makebox(0,0)[l]{$
H^0(\co(-C_j))\ \stackrel{\alpha_k}{\longrightarrow}
 \ H^0(\co_{\g_k}(-C_j))$}}
\put(0,25){\makebox(0,0)[l]{$
H^0(\co(-C_j))\ \stackrel{\alpha}{\longrightarrow}
  H^0(\co_{Z_{B_{j+1}}}(-C_j)) \to
H^1(\co(-C_{j+1}))\to $}}
\put(150,0){\makebox(0,0)[l]{$
 H^1(\co(-C_{j}))\to H^1(\co_{Z_{B_{j+1}}}(-C_j))\to 0 $}}
\put(30,35){\vector(0,1){30}}
\put(150,35){\vector(0,1){30}}
\put(40,50){\makebox(0,0){$=$}}
\put(160,50){\makebox(0,0){$\beta$}}
\end{picture}

By part ($d$): $h^1(\co(-C_t))=m-t$, hence 
 by dimension argument $h^1(\co_{Z_{B_{j+1}}}(-C_j))\not=0$. Since 
$\chi(\co_{Z_{B_{j+1}}}(-C_j))=0$ and $\beta$ is injective, we obtain that,
in fact, $\beta$ is an isomorphism, and $\alpha$ is
onto; hence $\alpha_k $ is onto as well.
Finally, the surjectivity of $\alpha_k$ for all $\g_k\subset B_{j+1}$ 
implies ($f$).

For $f\Rightarrow b$ notice that
the image of the elements 
$1,f_0,\ldots, f_{i-1}$ are linearly independent in  $H^0(\co_{C_i})$,
this together with (2.19) proves ($b$).

$f\Leftrightarrow g$. For any $\g_l\subset B_{j+1}$ consider the composed map
$$H^0(M,\co(-C_j))\stackrel{u}{\to}H^0(\co_{C'_{j+1}}(-C_j))
\stackrel{v}{\to}H^0(\co_{\g_l}(-C_j)).$$
By (2.2) $u$ is onto. Hence
$v\circ u$ is onto if and only if  $v$ is onto. Notice that
$h^0(\co_{\g_l}(-C_j))=1$ always (if $(f)$ is true then see the proof of $a
\Rightarrow f$; if $(g)$ is true, then $\co_{\g_l}(-C_j)$ 
is a trivial line bundle). 

$g\Leftrightarrow h$ is easy: 
$\co_{C'_{j+1}}(-C_j)=
\otimes_{i\leq j} \co_{C'_{i+1}}(-\zi)|_{C'_{j+1}}$, and
$\co_{C'_{j+1}}(-\zj)=
\co_{C'_{j+1}}(-C_j)\otimes \co_{C'_{j}}(-C_{j-1})|_{C'_{j+1}}^{-1}$.
\ \  \ $\Box$

\vspace{2mm}

\noindent 2.21.\ {\bf Corollary.}\ 
{\em Assume that \x is an elliptic Gorenstein singularity. 

a)\ Then $p_g=m+1$ if and only if 
the line bundles $\co_{C'_{j+1}}(-Z_{B_{j}})$ are trivial
for all $1\leq j\leq m-1$.

b)\ If $p_g=m+1$,  then 
 $h^1(M,\co(-l\zn))=h^1(\zk,\co_{\zk}(-l\zn))=p_g-1$ for any $l\geq 1$.}\\
{\em Proof} of\ $(a)$.
We have to prove that in the Gorenstein case
the triviality of $\co_{C'_1}(-\zbz)$ follows from the triviality of the 
other line bundles. Indeed, if all the other line bundles are trivial, then
$(2.20.a\Leftarrow h)$, applied for $(M/B_1,p_1)$,
gives  $h^1(\co_{C'_1})=m$. But the Gorenstein property implies that
for any positive cycle $D<\zk$ the strict inequality $h^1(\co_D)<h^1(\co_{\zk})
=p_g$ holds (see e.g. \cite{MR}, page 109).
Hence $h^1(\co_{\zk})>h^1(\co_{C'_1})=m$. Therefore $p_g=m+1$.
Now, we prove $(b)$. Using (2.2)
$h^1(\co(-l\zn))=h^1(\co_{\zk}(-l\zn))=
h^1(\co_{C'_1}(-l\zn))$
for $l\geq 1$.  But the line bundle
$\co_{C'_1}(-l\zn)$ is trivial by (2.20.g), hence 
its first Betti--number is $p_g-1$ by (2.20.e). \ \ \ $\Box$

\vspace{1mm}

\noindent 
In the next sections, (2.21.b) will replace some vanishing theorems
which were used in the classical case of  minimally elliptic singularities
(i.e. when  $p_g-1=0$). 

\vspace{1mm}

\noindent 2.22.\ {\bf Definition.} {\em In the sequel, we call the line bundles
$\co_{C'_{j+1}}(-\zj)$ ($0\leq j\leq m-1$) ``obstruction line bundles''. }

\vspace{1mm}

The next example shows that the property  $p_g=m+1$ is not
always true (even if we deal  with Gorenstein singularities).\\
2.23.\ {\bf Example.}\ 
Consider $(X_1,0)=\{x^2+y^3+z^{18}=0\}\subset (\bc^3,0)$ and 
$(X_2,0)=\{z^2=y(x^4+y^6)\}\subset (\bc^3,0)$ (cf. also \cite{Laufer87,Lauferbook} and \cite{Yau1}, page 291). Using the method of \cite{Lauferbook}, Chapter II,
(or the Appendix of \cite{Nemethi}) one can verify that the minimal resolution
graph in both cases is:

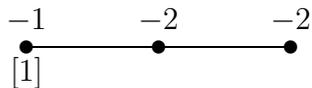
\begin{picture}(300,50)(0,0)
\put(100,25){\circle*{5}}
\put(150,25){\circle*{5}}
\put(200,25){\circle*{5}}
\put(100,25){\line(1,0){100}}
\put(100,35){\makebox(0,0){$-1$}}
\put(150,35){\makebox(0,0){$-2$}}
\put(200,35){\makebox(0,0){$-2$}}
\put(100,15){\makebox(0,0){$[1]$}}
\end{picture}

On the other hand, the geometric genus in the first case is $p_g=3$, and in the
second case is $p_g=2$ (this can be verified using e.g. \cite{Laufer77b} or 
\cite{Nemethi}). 

Denote by $\g_0,\g_1,\g_2$ the irreducible exceptional divisors (starting from left). Then $\g_0=E=\zbk$, $\g_0+\g_1=\zbe$, $\g_0+\g_1+\g_2=\zn$, and
$\zk=3\g_0+2\g_1+\g_2$ (hence $m=2$). Moreover $\chi(\zn)=0$, hence the
singularities are elliptic.
It is very instructive to check property (2.20.f) for these examples.
For a function $f:(X,0)\to (\bc,0)$ we denote its divisor supported by
$\g$ by $(f)_{\g}:=\sum_im_{\g_i}(f\circ \pi)\g_i$. 

In the first case, $\zn=(z)_{\g}$. 
This shows that in (2.20.f) (a possible choice is) $f_0=z$ and $f_1=z^2$.
In the second case there is no function $f_0$ with $m_{\g_0}(f_0)=
m_{\g_0}(\zn)=1$, hence (2.20.f) fails (in particular, there is no function
$f$ with $(f)_{\g}=\zn$.)
\\

\noindent {\bf 3.\ Characterization of the trivial line bundles in $Pic(C'_j)$.}\\

As usual, for any positive cycle $D$, we denote the isomorphism classes
of invertible sheaves over $\co_D$ by $Pic(D)$ ($=H^1(\co^*_D)$). The 
kernel of the degree
map $deg:Pic(D)\to \bz^{\#\{\g_i|\g_i\subset |D|\}}$ defined by $\cl\to
\{deg_{\g_i}\cl\}_i$ is denoted by $Pic^0(D)$. 

If $D=E$ is a minimally elliptic cycle, then
$$(3.1)\ \ \left\{\begin{array}{l}
\cl\in Pic^0(E)\ \mbox{is trivial}\Leftrightarrow H^0(E,\cl)\not=0;\\
\mbox{moreover, if $\cl$ is trivial then $h^0(E,\cl)=1$.}\end{array}\right. $$
This basically follows from the 2--connectivity of $E$ (see e.g. \cite{MR}
page 82). But, in general, it is not really easy to provide similar result
for $Pic^0(D)$. Fortunately, in the case of an elliptic singularity,
for any $0\leq j\leq m$, $Pic^0(\zj)$ is as simple as $Pic^0(E)$. Indeed, 
consider the ``exponential cohomology sequences'' (see, e.g. \cite{BPV}
page 49) of the cycles $\zj$ and $E$.
Then  the natural maps $H^1(|\zj|,\bz)\to H^1(|E|,\bz)$ and
$H^1(\co_{\zj})\to H^1(\co_E)$ are isomorphism (the first from 2.9, or by a
similar argument; for the second use 2.16 for $D=0$). Hence 
$Pic^0(\zj)\to Pic^0(E)$ is an isomorphism as well.
Moreover, if we consider a computation sequence $\{Z_l\}_l$
which connects $E=\zm$ and $\zj$
(cf. 2.15), then from the exact sequences
$0\to \co_{\g_{i_l}}(-Z_l)\otimes \cl\to\cl|Z_{l+1}\to \cl|Z_l\to 0$
follows that $H^0(\zj,\cl)=H^0(E,\cl|_E)$. Therefore, for any line bundle
$\cl\in Pic^0(\zj)$:
$$(3.2)\ \ \left\{\begin{array}{l}
h^0(\zj,\cl)\leq 1,\ \mbox{and}\\
h^0(\zj,\cl)=1\ \mbox{if and only if $\cl$ is trivial}.\end{array}\right. $$
3.3.\ {\bf Theorem.}\ {\em Assume that the elliptic numerical Gorenstein
singularity \x
satisfies $p_g=m+1$; and fix an integer $0\leq j\leq m$. Then the followings
are equivalent.

a)\ $\cl\in Pic^0(C'_j)$ is trivial;

b)\ $h^0(C'_j,\cl)\geq m-j+1.$

\noindent 
Actually, if $\cl$ is trivial, then $h^0(C'_j,\cl)$ is exactly $m-j+1$.}\\
{\em Proof.}\ $a\Rightarrow b$ and the equality follow from 
(2.20; $a\Rightarrow e$) and (2.14). We prove 
$b\Rightarrow a$ by descending induction over $j$. If $j=m$, then $C'_j=E$,
hence (a) follows from (3.1). Assume that $b\Rightarrow a$ is true for $j+1$,
and $h^0(C'_j,\cl)\geq m-j+1$. Consider the exact sequence
$$0\to \co_{C'_{j+1}}(-\zj)\otimes \cl\to \cl\to\cl|\zj\to 0.$$
Then $h^0(\cl|\zj)\leq 1$ by (3.2), hence $h^0:=h^0(\co_{C'_{j+1}}
(-\zj)\otimes \cl)\geq m-j$. Hence, by the inductive step $\co_{C'_{j+1}}
(-\zj)\otimes\cl$ is trivial in $Pic^0(C'_{j+1})$, and $h^0=m-j$. Therefore
$h^0(\cl|\zj)$ must be 1, $\cl|\zj$ must be trivial, and $r:H^0(\cl)\to H^0(
\cl|\zj)$ onto. In particular, there exists $s\in H^0(\cl)$ such that 
$r(s)$ generates $\co(\cl|\zj)$. Since $|\zj|=|C'_j|$, it follows that $s$
generates the sections of $\cl$ as well.\ \ \ $\Box$\\
3.4.\ {\bf Corollary.}\ {\em Assume that \x is an elliptic Gorenstein
singularity, and we fix an integer $0\leq j\leq m$. Assume that 
$h^1(\co_{C'_j})=m-j+1$. Then 
$\co_{C'_j}(-C_{j-1})$ is trivial in $Pic^0(C'_j)$.

For example, for $m=j$ we obtain that
$\co_E(-C_{m-1})$ is trivial in $Pic^0(E)$ (the assumption 
$h^1(\co_E)=1 $ is automatically satisfied). For 
$j=m-1$, if contracting $B_{m-1}$ one has $p_g((M/B_{m-1},p_{m-1}))=2$, then 
$\co_{C'_{m-1}}(-C_{m-2})$ is trivial in $Pic^0(C'_{m-1})$.}\\
{\em Proof.}\ First notice that contracting $B_j$ we obtain a numerical
 Gorenstein singularity with canonical cycle $C'_j$ (cf. 2.12), and with the
``extremal property'' $p_g(M/B_j,p_j)=
$ length of the elliptic sequence of $B_j=m-j+1$. Therefore, (3.3) for this singularity gives:
$\cl\in Pic^0(C'_j)\ \mbox{is trivial}\Leftrightarrow h^0(C'_j,\cl)=m-j+1.$
But by Serre duality
$h^0(\co_{C'_j}(-C_{j-1})=h^1(\co_{C'_j})=m-j+1$.\ \ \ \ $\Box$

\vspace{2mm}

\noindent 
3.5.\ {\bf Second characterization of $p_g=m+1$.}\ {\em Assume that \x
is a numerical Gorenstein elliptic singularity. Consider its minimal 
resolution $M\to X$ and the singularity $(M/B_j,p_j)$ 
obtained by the contraction of $B_j$ ($0\leq j\leq m$). 
Then \x satisfies the extremal property 
$p_g=m+1$ if and only if the singularities $(M/B_j,p_j)$ are Gorenstein
for all $0\leq j\leq m-1$.
(Notice that $(M/B_m,p_m)$ is automatically Gorenstein by 2.7.)}\\
{\bf Proof.}\ ``$\Rightarrow$''\ 
First we prove that a numerical Gorenstein singularity \x with $p_g=m+1$
is Gorenstein. (This fact was proved by \y  in \cite[(3.11)]{Yau1}; we present
here a very short proof in the spirit of \cite[(4.21)]{MR}; cf. also with
the papers of L. B\u adescu \cite{B1,B2}.)

Consider the line bundle $\cl:=\co_{\zk}(\zk+K)\in Pic^0(\zk)$. By Serre
duality $h^0(\cl)=h^1(\co_{\zk})$, hence $h^0(\cl)=p_g=m+1$. Therefore, by (3.3), $\cl$
is trivial. Moreover, by (2.2) $h^1(\co(K))=0$, hence
 $H^0(M,\co(\zk+K))\to H^0(\cl)\approx
H^0(\co_{\zk})$ is onto. Therefore, there is a global section of
$\co(\zk+K)$ which has no zeros in the neighbourhood of $\g$, hence
$\zk+K$ is linearly equivalent to the zero cycle. But this is one of the 
characterizations of the Gorenstein property.

Now, $M/B_j$ is an elliptic
numarical Gorenstein singularity
with numerical cycle $C'_j$ and the length 
of its elliptic sequence$=h^1(\co_{C'_j})=m-j+1$ (cf. 2.11 and 2.20.e). 
Therefore, the above fact
applied for the singularities $M/B_j$ ends the proof.

``$\Leftarrow$''  By (2.20.h),
it is enough to prove that the obstruction line bundles $\cl_j:=\co_{C'_j}
(-Z_{B_{j-1}})$ are trivial for $1\leq j\leq m$. We will prove this by 
induction. Assume that for a fixed $1\leq j_0\leq m$ the line bundles
$\cl_j$ are trivial for $j>j_0$
(if $j_0=m$ then this assumption is vacuous). Since the singularity
$M/B_{j_0-1}$ is Gorenstein, by (2.20.h) and (2.21.a) (applied for
$M/B_{j_0-1}$), the triviality of 
$\{\cl_j\}_{j>j_0}$ implies the triviality of $\cl_{j_0}$.\ \ \ $\Box$\\

\noindent 
3.6.\ {\bf Examples.}\ Assume that \x is a numerical Gorenstein elliptic
singularity.

(a)\ If $m=1$ then $p_g\leq 2$, and:
$$p_g=2\Leftrightarrow \co_E(-\zn) \ \mbox{is trivial}\ \Leftrightarrow
\mbox{\x is Gorenstein.}$$
Examples with $p_g=1$ exists (actually this is the generic case!).
E.g. take the minimal resolution of
any Gorenstein singularity with $p_g=2$ and $m=1$. Then deform
its analytic structure.  Generically we obtain a non--Gorenstein singularity 
with $p_g=1$ (cf. \cite{Laufer77}, page 1279).

(b)\ Assume that $m=2$. Then in the numerical Gorenstein case we have two 
obstruction line bundles. For the generic analytic structure $p_g=1$ as
above. If \x is Gorenstein, then $p_g=h^1(\co_{\zk})>h^1(\co_E)=1$, hence
$p_g\geq 2$. From (2.21.a), $p_g=3$ if and only if $\co_E(-\zbe)$ is trivial.
From (3.4) $\co_E(-\zn-\zbe)$ is trivial, hence {\bf if \x is Gorenstein, then:}
$$ p_g=3\Leftrightarrow \ \co_E(-\zn)\ \mbox{is trivial}\Leftrightarrow
(M/B_1,p_1) \ \mbox{is Gorenstein}.$$\\

\noindent {\bf 4.\ Torsion properties of the ``obstruction line bundles''.}\\

\noindent The main result of this section is the following:\\
4.1.\ {\bf Theorem.}\ {\em Consider the minimal resolution of
an elliptic Gorenstein singularity \x
with elliptic sequence $\{\zj\}_{j=0}^m$ ($m\geq 1$), and fix an
integer $1\leq k\leq m$.

Assume that for $k+1\leq j\leq m$ the ``obstruction line bundles''
$\co_{C'_j}(-Z_{B_{j-1}})\in Pic^0(\co_{C'_j})$ are trivial. (If $k=m$
then this assumption is vacuous.) Then $\co_{C'_k}(-Z_{B_{k-1}})\in Pic^0(
\co_{C'_k})$ has finite order.

Actually, there exist integers $1\leq l_r\leq k-r$ \ ($0\leq r\leq k-1$)
such that for any $0\leq r\leq k-1$ the line bundle $\co_{C'_k}(l_0\cdots
l_r\ Z_{B_r})$ is trivial. In particular, the order of 
$\co_{C'_k}(-Z_{B_{k-1}})$ in $ Pic^0(\co_{C'_k})$ is not greater than $k!$.}\\

\noindent {\em Proof.}\ First assume that $k=1$. 
Then $\co_{C'_1}(-\zbz)$ is trivial by (2.21.a) and $(2.20.a\Leftrightarrow h)$.
So, in the sequel we will assume that $k\geq 2$.

Since  $(M/B_k,p_k)$ is a numerical
Gorenstein singularity with canonical class $C'_k$ (cf. 2.11.d),
the assumption together with (2.20) (applied for this singularity) provide
that:
$$(4.2)\hspace{4.5cm}
h^1(\co_{C'_k})=m-k+1.\hspace{4.5cm}$$
\noindent 
4.3.\ {\bf Lemma.} \ {\em For any $0\leq r\leq t\leq k-2$ one has
$h^1(\co(-C_t-Z_{B_r})) \geq m-k+1$.}\\
{\em Proof.}\ Using (2.2), Serre duality, (2.14), and in the last step
$C'_{t+2}\geq C'_k$, one has:
$h^1(\co(-C_{t+1}))\stackrel{2.2}{=}h^1(\co_{C'_{t+2}}(-C_{t+1}))
\stackrel{SD}{=}h^0(\co_{C'_{t+2}})\stackrel{2.14}{=}h^1(\co_{C'_{t+2}})
\geq h^1(\co_{C'_k})$. Hence, by (4.2):
$$(4.4)\hspace{4cm}
h^1(\co(-C_{t+1})\geq m-k+1.\hspace{4cm}$$
Now, consider a computation sequence $\{Z_l\}_l$ which connects $Z_{B_{t+1}}$
with $Z_{B_r}$ (cf. 2.15) (i.e. $Z_{l+1}=Z_l+\g_{i_l}$ where $\g_{i_l}$
is  smooth rational curve with $\g_{i_l}Z_l=1$). In the exact sequence
$0\to \co(-C_t-Z_{l+1})\to \co(-C_t-Z_l)\to \co_{\g_{i_l}}(-C_t-Z_l)\to 0$
the Chern number $\g_{i_l}(-C_t-Z_l)$ is $-\g_{i_l}C_t-1\geq -1$ (cf. 2.11.e).
Hence $H^1(\co(-C_t-Z_{l+1}))\to H^1(\co(-C_t-Z_{l}))$ is onto for any $l$,
in particular, $h^1(\co(-C_t-Z_{B_r}))\geq h^1(\co(-C_t-Z_{B_{t+1}}))$. The lemma follows from this and (4.4). \ \ \ $\Box$\\
\noindent Now consider the exact sequence:
$$(4.5)\hspace{1cm}
0\to \co(-\zk-Z_{B_r})\to \co(-C_t-Z_{B_r})\to \co_{C'_{t+1}}(-C_t-Z_{B_r})\to 0.$$
4.6.\ {\bf Lemma.}\ $H^1(\co(-\zk-Z_{B_r})=0$.\\
{\em Proof.}\ Consider a computation sequence $\{Z_l\}_l$ which connects 
$Z_{B_r}$ and $Z_{B_0}=\zn$ (cf. 2.15). Consider the exact sequences:
$$0\to \co(-\zk-Z_{l+1})\to \co(-\zk-Z_l)\to \co_{\g_{i_l}}(-\zk-Z_l)\to 0$$
and notice that $\g_{i_l}(-\zk-Z_l)=-\g_{i_l}\zk-1\geq -1$ (cf. 2.1).
Hence $H^1(\co(-\zk-\zn))\to H^1(\co(-\zk-Z_{B_r}))$ is onto.
But $H^1(\co(-\zk-\zn))=0$ (by 2.2).\ \ \ $\Box$\\
Now, from (4.5) and duality:
$$H^1(\co(-C_t-Z_{B_r}))=H^1(\co_{C'_{t+1}}(-C_t-Z_{B_r}))=
H^0(\co_{C'_{t+1}}(C'_{t+1}+C_t+Z_{B_r}-\zk)).$$
Hence (4.3) reads as follows: for $0\leq r\leq t \leq k-2$:
$$h^0(\co_{C'_{t+1}}(C'_{t+1}+C_t+Z_{B_r}-\zk))\geq m-k+1.$$

Now, for any positive cycle $0<D\leq C'_{t+1}$ write $o(D):=\co_D(D+C_t+
Z_{B_r}-\zk)$. For any irreducible exceptional divisor $\g_i\subset
|D|$, one has the exact sequence:
$$0\to o(D-\g_i)\to o(D)\to\co_{\g_i}(D+C_t+Z_{B_r}-\zk)\to 0.$$
The Chern number $\g_i(D+C_t+Z_{B_r}-\zk)=\g_i(D-\zk)$, for 
$\g_iZ_{B_r}=\g_iC_t=0$. Therefore, if 
$\g_i(D-\zk)<0$, then 
$H^0(o(D-\g_i))\to H^0(o(D))$ is an isomorphism. 
In the sequel, we will construct a decreasing sequence of cycles as follows.
The starting cycle is $D=C'_{t+1}$. If $D$ is already constructed, and
there exists $\g_i\subset |D|$ with $\g_i(D-\zk)<0$, then the next term is
$D-\g_i$. During this step $h^0(o(D))$ will stay constant.

If at a moment we reach a cycle $D>0$ with the property $\g_i(D-\zk)\geq 0$
for all $\g_i\subset |D|$, then by (2.13) in fact
$\g_i(D-\zk)=0$ for all 
$\g_i\subset |D|$ and 
$D$ must be one of the cycles $C'_s$. 
If there is a $\g_i\subset |D|$ such that
$$r_i:H^0(\co_{D}(D+C_t+Z_{B_r}-\zk))\to H^0(\co_{\g_i}
(D+C_t+Z_{B_r}-\zk))$$
is the trivial map, then the term following $D$ is $D-\g_i$ for this $\g_i$. 
Hence, in this case again,
$h^0(o(D))=h^0(o(D-\g_i))$. If the above map $r_i$ is nontrivial
 for all $\g_i\subset
|D|$, then  the term following $D$ is $D-\g_i$, where 
$\g_i$ arbitrary with $\g_i\subset |D|$. 
In this later case $r_i$ has rank one, hence $h^0(o(D))$ will decrease by one.
If $D=0$ we stop.

Since finally $h^0(o(D))$ must decrease to zero, and $h^0(o(C'_{t+1}))\geq 
m-k+1$, $h^0(o(D))$ must decrease at least $m-k+1$ times. Therefore, the first time when we decrease $h^0$ must happen for a cycle $D=C'_s$ with
$s\leq k$ (i.e. the length of the sequence $\{C'_s,\ldots,C'_m\}$
must be at least $m-k+1$). Obviously $C'_s\leq C'_{t+1}$
(because we deal with cycles $D\leq C'_{t+1}$), hence: $s\geq t+1$. 

This shows that there exists $t+1\leq s\leq k$ such that 
$$H^0(\co_{C'_s}(C'_s+C_t+Z_{B_r}-\zk))\to H^0(\co_{\g_i}
(C'_s+C_t+Z_{B_r}-\zk))=\bc$$
is onto for all $\g_i\subset |C'_s|$.  In other words, $
\co_{C'_s}(C'_s+C_t+Z_{B_r}-\zk)$ is a trivial line bundle. Taking its 
restriction to $C'_k$, we obtain:
$$(4.7)\hspace{1cm}
\left\{\begin{array}{l}
\mbox{for any pair $(r,t)$ with $0\leq r\leq t\leq k-2$,}\\
\mbox{there exists $s$ with $t+1\leq s\leq k$, such that}\\
\co_{C'_k}(C'_s+C_t+Z_{B_r}-\zk)\ \mbox{is a trivial line bundle
in}\ Pic^0(C'_k).
\end{array}\right. $$
4.8.\ {\bf Lemma.}\ {\em For any $0\leq r\leq k-1$, there exists $l_r$ with
$1\leq l_r\leq k-r$, such that }
$$\co_{C'_k}(l_0\cdots l_r\ Z_{B_r})\ \mbox{is trivial in}\ Pic^0(C'_k).$$
{\em Proof.}\ First we prove for $r=0$. If there is at least one $s$ in (4.7)
with $s=t+1$, then 
$C'_s+C_t-\zk=0$, hence $\co_{C'_k}(Z_{B_0})$ is trivial. Now, 
assume that  $s\geq t+2$ for any $t$.  This means that for any 
$0\leq t\leq k-2$, there exists $t+2\leq s(t)\leq k$ such that:
$$(*_t)\hspace{2cm}
\co_{C'_k}(\zbz)\approx\co(Z_{B_{t+1}}+\cdots+ Z_{B_{s(t)-1}}).$$
Now, we set $t_1:=0$. If $s(t_1)=k$ we stop; otherwise consider
$t_2:=s(t_1)-1$. If $s(t_2)=k$ we stop, otherwise we continue
until we obtain $s(t_u)=k$ for some $u\geq 1$. If we multiply
 the identities $(*_t)$ for $t_1,\ldots, t_{u}$,
we obtain:
$$\co_{C'_k}(u\zbz)\approx\co_{C'_k}(\zbe+\cdots +Z_{B_{k-1}}).$$
Here $u$ satisfies $1\leq u\leq k-1$ ($u=k-1$ exactly when $s(t)=t+2$ for 
every $t$,
i.e. when  we have to consider all the isomorphisms $(*_t)$).
But from (3.4) and (4.2):
$$(4.9)\hspace{4cm}
\co_{C'_k}(\zbz+\cdots +Z_{B_{k-1}})\ \mbox{is trivial},\hspace{4cm}$$
hence  with the notation $l_0:=u+1$ one has 
$\co_{C'_t}(l_0\zbz)$ is trivial.\\

Now, take $r=1$ and consider again (4.7). 
If $s=t+1$ for at least one $t$, then $\co_{C'_k}(\zbe)$ is trivial.
If $s\geq t+2$ for all $1\leq t\leq k-2$, by a similar argument as above
one obtains that 
$\co_{C'_k}(u\zbe)\approx \co_{C'_t}(\zbk+\cdots+Z_{B_{k-1}})$
for some $1\leq u\leq k-2$ ($u=k-2$ occurs exactly when we have to consider 
all the identities corresponding to $1\leq t\leq k-2$). This together with 
(4.9) gives that
$\co_{C'_k}(\zbz+(u+1)\zbe)\ \mbox{is trivial}.$
Now consider its $l_0^{th}$ power and take $l_1:=u+1$.

Similar argument can be used for all $0\leq r\leq k-2$. For $r=k-1$, just take 
(4.9) again and take its $l_0\cdots l_{k-2}$-th power.\ \ \ $\Box$

\vspace{2mm}

\noindent 4.10. {\bf Example.}  (cf. 3.6.b)\ Assume that \x is elliptic
Gorenstein, and $m=2$. Recall that  $\co_E(-\zbe-\zn)\approx\co_E$. 
The above theorem (for $m=k=2$)
says that the order of $\co_E(-\zn)$ divides 2, hence:
$$\co_E(-\zn)\ \mbox{is trivial}\Leftrightarrow p_g=3;\hspace{1cm}
\co_E(-\zn)\ \mbox{has order 2}\Leftrightarrow p_g=2.$$
One of the main results of the present paper is the following:

\vspace{2mm}

\noindent 
4.11.\ {\bf Theorem.}\ {\em Assume that \x is an elliptic Gorenstein
singularity with  $H^1(\g,\bz)=0$. Then $p_g=m+1=$ the length of the elliptic 
sequence in  the minimal resolution of $(X,p)$. 
In particular $p_g$ is a topological invariant.}\\
{\em Proof.}\ If $H^1(\g,\bz)=0$, then by the exponential exact sequence
(\cite{BPV}, page 49) $Pic^0(D)\approx H^1(\co_D)$ for any positive cycle $D$,
in particular, it is torsion free. Therefore, (4.1) (as an inductive step)
proves that all the obstruction line bundles are trivial. Hence $p_g=m+1$ 
by (2.20).\ \ \ $\Box$\

\vspace{2mm}

\noindent 4.12.\ {\bf Remark.}\ If $m\geq 1$, then in (4.11) both assumptions 
are necessary. Without Gorenstein condition $p_g$ generically is one
(\cite[page 1279]{Laufer77}). Moreover,
Example (2.23) shows that $p_g<m+1$ can occur if $H^1(\g,\bz)\not=0$.\\
4.13.\ {\bf Corollary.}\ {\em For  a Gorenstein  singularity \x with
$H^1(\g,\bz)=0$ the  following facts are equivalent:}
$$p_g=2 \ \Leftrightarrow \ \chi(\zn)=0 \ \mbox{and} \ \zk=\zn+E.$$
{\em (Notice that the right hand side is completely topological.)}\\
{\em Proof.}\ Indeed, if \x is Gorenstein with $p_g=2$, then $h^1(\zn)$ cannot 
be zero (because of 2.4), cannot be greater than or equal to
two (that whould imply that
$\zn=\zk$, which characterise the minimally elliptic singularities, cf. 2.7),
hence it is one. Therefore \x is elliptic.
The rest follows from the above results.\ \ \ $\Box$\\
\noindent 4.14. \
{\bf Remark.}\ \y in \cite{Yau1}, Theorem B, proved that a Gorenstein  
singularity with $p_g=2$ is elliptic. Moreover, in \cite{Yau3}, he proved
also that  a {\em hypersurface} singularity with $H^1(\g,\bz)=0$ and $p_g=2$
 satisfies $m=1$. His proof (of this second fact)
is based on the classification of all possible 
dual resolution graphs of hypersurface singularities with $p_g=2$
(250 cases).

\vspace{2mm}

\noindent Notice that  the ``second characterization'' (3.5) gives:\\
4.15.\ {\bf Corollary.}\ {\em 
If \x is an elliptic Gorenstein
singularity with  $H^1(\g,\bz)=0$ then 
the singularities $(M/B_j,p_j)$ are Gorenstein for all $0\leq j\leq m$.}\\

\noindent {\bf 5. \ The multiplicity.}\\

\noindent
5.1.\ Assume that \x is an elliptic Gorenstein singularity with $p_g=m+1$.
In this section we will prove that its multiplicity $mult(X,p)$ is a 
topological invariant, in fact depends only on $\zn^2$.

First notice that by (2.20.f) there exists $
f_0\in H^0(M,\co(-\zn))$ such that for any $\g_l\subset B_1$,
the order of vanishing of $f_0$ on $\g_l$ is exactly $m_{\g_l}(\zn)$.
In particular:
$$(5.2)\hspace{1cm}
H^0(M,\co(-\zn))\to H^0(\co_{\g_l}(-\zn))=\bc\ \
\mbox{is onto for all $\g_l\subset B_1$.}$$
By (2.2), if $|\co_{\zk}(-\zn)|$ is free, then $|\co(-\zn)|$ is free too. 
By \cite[page 426]{Wagreich}, 
if $|\co(-\zn)|$ is free, then $m_p\co_M=\co(-\zn)$, in which case 
$mult(X,p)=-\zn^2$. In particular:
$$(5.3) \hspace{5mm} \zn^2=-1\ \Rightarrow 
-\zn^2\not= mult(X,p)
\Rightarrow 
\ |\co_{\zk}(-\zn)|\ \mbox{is not free.}$$
If $|\co(-\zn)|$ is not free, then $mult(X,p)>-\zn^2$. E.g. if 
$m_p\co_M=m_Q\co(-\zn)$, where $m_Q$ is the maximal ideal of a smooth
point $Q$ of $\g$, then $mult(X,p)=-\zn^2+1$. 

The next theorem generalizes Laufer's result about the multiplicity of 
minimally elliptic singularities \cite{Laufer77}.

\vspace{2mm}

\noindent 5.4.\ {\bf Theorem.}\  {\em Assume that \x is an elliptic 
Gorenstein singularity with $p_g=m+1$. 
(If $m=0$ or $H^1(\g,\bz)=0$ then the last assumption is satisfied.)
Then:

If $\zn^2\leq -2$, then $m_p\co_M=\co(-\zn)$, hence $mult(X,p)=-\zn^2$;

If $\zn^2=-1$, then $m_p\co_M=m_Q\co(-\zn)$ for some smooth point $Q$
of $\g$, and $mult(X,p)=2$.}

\vspace{2mm}

\noindent
We prove theorem (5.4) in several steps. First we prove the converse of (5.3).
For simplicity, in the sequel we will write $\cl:= \co_{\zk}(-\zn)$.

\vspace{2mm}

\noindent 5.5.\ {\bf Proposition.} \ {\em Assume that $|\cl|$
is not free. Then $\zn^2=-1$.}

\vspace{2mm}

\noindent
{\em Proof.}\ By (5.2), all the basepoints must lie outside of $B_1$. 
Since all the connected components of $\overline{\g\setminus B_1}$ are trees 
whose vertices correspond to smooth rational curves (cf. 2.8),
we can fix a basepoint $Q$ such that $Q$ is ``as close to $B_1$ as possible''.
 More precisely, we will fix a basepoint $Q$ such that
there exist irreducible exceptional divisors $\{\g'_i\}_{i=0}^l$ with 
$$(5.6)\hspace{1cm}\left\{\begin{array}{l}
\mbox{$\g'_i\not\subset B_1$ and  $\g'_i$ is smooth rational curve
$(0\leq i\leq l)$};\\
\mbox{$Q\in \g'_0,\ \g'_i\g'_{i+1}=1$ for $0\leq i\leq l-1$, and 
$\g'_lB_1=1$;}\\
\mbox{ there is no basepoint of $|\cl|$ on $\cup_{i>0}\g'_i$.}
\end{array}\right.$$
With this choice of $Q$, $\g'_0\zn$ cannot be zero. Indeed, since $Q$ is a 
basepoint, the evaluation map
$H^0(\cl)\stackrel{\alpha}{\longrightarrow}H^0(\co_{\g'_0}(-\zn))
\stackrel{e}{\longrightarrow}\bc_Q$
is trivial. If, by assumption, $\g'_0\zn=0$ then
 $e$ must be an isomorphism, hence
$\alpha$ must be trivial  as well. This shows that $\g'_0\cap\g'_1$ is also
a basepoint, which contradicts the above choice of $Q$ in (5.6). For later
reference:
$$(5.7)\hspace{5cm}\g'_0\zn<0.\hspace{3cm}\hspace{5cm}$$
Since the evaluation map $H^0(\zk,\cl)\to\bc_Q$ is zero, from the exact 
sequence $0\to m_Q\cl\to\cl\to\bc_Q\to 0$ and (2.21.b),
one gets $h^1(m_Q\cl)=h^1(\cl)+1=p_g$.

In the next paragraph, we follow the beginning of the proof of (4.23)
of \cite{MR}. 
Take $0<D\leq \zk$ minimal with the property $h^1(D,m_Q\cl|_D)=p_g$. Minimality means:
$$h^1(D-\g_i,m_Q\cl|_{D-\g_i})\leq p_g-1\ \mbox{for any}\ \g_i\subset |D|.$$
Let $K_i$ be the kernel of $m_Q\cl|_D\to m_Q\cl|_{D-\g_i}$. 
Outside $Q$, clearly
$K_i=\cl\otimes \co_{\g_i}(-D+\g_i)$. Actually, if $Q\in|D-\g_i|$, then even at $Q$ one has
$K_i=\cl\otimes \co_{\g_i}(-D+\g_i)$. If $Q\in|D|$ but $Q\not\in|D-\g_i|$
(i.e. if $Q\in\g_i\setminus\{\mbox{other components of $D$}\}$ and $m_{\g_i}(D)
=1$), then 
$K_i=\cl\otimes m_Q\co_{\g_i}(-D+\g_i)$.
We define $\delta_i:=0$ in the first two cases, and $\delta_i:=1$ in the third case.

Now, the minimality of $D$ implies that for all $i$: $h^1(D,K_i)\not=0$,
in other words:
$$(5.8)\hspace{3cm}
\g_i\cl+\g_i(\zk-D)\leq \delta_i\ \mbox{for all $\g_i\subset |D|$}.\hspace{2cm}$$
5.9.\ {\bf Lemma.}\ {\em The situation when 
$\g_i\cl + \g_i(\zk-D)=0$ for all $\g_i\subset |D|$ cannot occur.}\\
{\em Proof.}\ Assume the contrary.
Since $\g_i\cl\geq 0$ for all $\g_i$, and $\chi(D)\geq 0$, the relation
$$0\leq D(\zk-D)=\sum m_{\g_i}(D)\g_i(\zk-D)=-\sum m_{\g_i}(D)\g_i\cl \leq 0$$
says that $\g_i\cl=\g_i(\zk-D)=0$ for all $\g_i\subset |D|$. In particular,
$D=C'_s$ for some $0\leq s\leq m$ (by 2.13), and $D\cl=C'_s(-\zn)=0$. Since
$C'_0(-\zn)=-\zn^2\not=0$, we must have $s\geq 1$, i.e. $|D|\subset B_1$.
On the other hand, (5.2) guarantees that $Q\not\in B_1$, hence $Q\not\in |D|$.
Therefore, $m_Q\cl|_D=\co_{C'_s}(-\zn)$ for some $s\geq 1$. Using (3.3) we obtain 
$h^1(\co_{C'_s}(-\zn))=h^0(\co_{C'_s}(-\zn))\leq m-s+1\leq p_g-1$, which
contradicts  $h^1(D,m_Q\cl|_D)=p_g$.\ \ \ \ $\Box$

\vspace{2mm}

\noindent 
The above lemma implies that there exists $\g_{i_0}\subset |D|$ with 
$\g_{i_0}\cl+\g_{i_0}(\zk-D)=\delta_{i_0}=1$. (In particular, either
 $Q=\g_{i_0}\cap\g'_0 $ and $\g'_0\not\subset |D|$,
or $\g_{i_0}=\g'_0$, cf. 5.6. The first case will be eliminated later.)
Summing (5.8) over all components of $D$ (recall $m_{\g_{i_0}}(D)=1$)
we obtain 
$D\cl+D(\zk-D)=1$. Since $D\cl\geq 0$ (by the definition of $\zn$)
and $D(\zk-D)$ is an {\em even} 
non--negative integer (by Riemann--Roch), we get:
$$(5.10)\hspace{4cm}
D\cl=1\ \mbox{and}\ D(\zk-D)=0.\hspace{4cm}$$
If we write $D=\sum D_i$, where $\{|D_i|\}_i$ are the the connected components
of $|D|$, since $\chi(D)$=0 and $\chi(D_i)\geq 0$ for all $i$ (cf. 2.5),
we must have $\chi(D_i)$ for all $i$. This shows that no $D_i$ lies in 
$\overline{\g\setminus B_1}$ (for all the components of $\overline{\g\setminus
B_1}$ support rational singularities whose positive cycles
 have $\chi>0$, cf. 2.8 
and 2.4). Therefore, all the exceptional divisors $\{\g'_i\}_{i=0}^l$ defined in
(5.6) must be components of $|D|$ and $\g'_0=\g_{i_0}$. Then (5.7) reads as
$$(5.11)\hspace{5cm}
\g_{i_0}\zn<0.\hspace{5cm}$$
Now, $\g_j\cl\geq 0$ for all $\g_j\subset |D|$ and $D\cl=1$ (cf. 5.10), hence
there is exactly one index $j_0$ with $\g_{j_0}\zn\leq -1$. (5.11)
guarantees that this $j_0$ is exactly $i_0$. Therefore, we get:
$$\begin{array}{lll}
\g_{i_0}\zn=-1, & \g_{i_0}(\zk-D)=0; & \\
\g_{j}\zn=0, & \g_{j}(\zk-D)=0 & \mbox{for $j\not=i_0;\ \g_j\subset |D|.$}
\end{array}$$
Again by (2.13): $D=C'_s$ for some $0\leq s\leq m$. Since $Q\in |D|$ but 
$Q\not\in  B_1$, the only possibility is $D=\zk$. This means
that $1=D\cl=\zk(-\zn)=-\zn^2$, which ends the proof of
(5.5).  \ \ \  $\Box$\\

\noindent 5.12.\ {\bf The special case $\zn^2=-1$.}\ \ 
All the possible dual (minimal) resolution graphs of elliptic numerical 
Gorenstein singularities with $\zn^2=-1$ are classified. This is done by
\y (in \cite{Yau2}, (2.3) and Table 1) based on Laufer's classification 
of the minimal resolution
graphs of the minimally elliptic singularities with $E^2=-1$ (\cite{Laufer77},
Table 1.) \\

\noindent 5.13.\ {\bf Proposition.} \ \cite{Yau2,Laufer77}\ {\em 
The minimal resolution graph of a numerical Gorenstein elliptic singularity
with $\zn^2=-1$ has the form $\g=B_1\cup \g_0$, where $\g_0$ is a smooth
rational curve with self--intersection $-2$, and it is unique
with  $\g_0\zn=-1$. Moreover, $m_{\g_0}(\zn)=m_{\g_0}(\zk)=1$. \ 
$B_1$ supports a singularity with  $\zbe^2=-1$. 
The curve $\g_0$
intersects only one irreducible exceptional divisor $\tilde{\g}$ of $B_1$,
and it intersects this exactly in one point. 
The curve $\tilde{\g}$ is the unique component of $B_1$ with 
the property $\tilde{\g}\zbe=-1$.

Therefore, by induction, $\g=|E|\cup\g_{m-1}\cup\cdots\cup\g_0$, 
where all curves $\{\g_i\}_{i=0}^{m-1}$ are smooth rational curves with 
self--intersection $-2$, $E^2=-1$, $\zj=E+\sum_{i\geq j}\g_i$, and the dual
graph has the following form:}

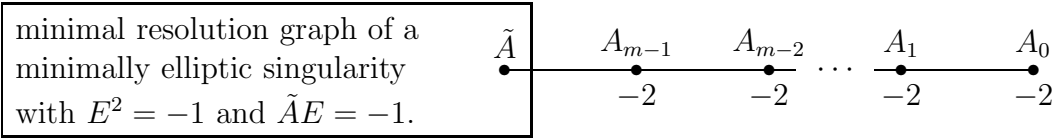
\begin{picture}(400,70)(30,-10)
\put(10,0){\framebox(200,50){}}
\put(15,40){\makebox(0,0)[l]{minimal resolution graph of a}}
\put(15,25){\makebox(0,0)[l]{minimally elliptic singularity}}
\put(15,10){\makebox(0,0)[l]{with $E^2=-1$ and $\tilde{\g}E=-1$.}}
\put(200,25){\circle*{4}}
\put(250,25){\circle*{4}}
\put(300,25){\circle*{4}}
\put(350,25){\circle*{4}}
\put(400,25){\circle*{4}}
\put(200,25){\line(1,0){110}}
\put(400,25){\line(-1,0){60}}
\put(200,35){\makebox(0,0){$\tilde{\g}$}}
\put(250,35){\makebox(0,0){$\g_{m-1}$}}
\put(300,35){\makebox(0,0){$\g_{m-2}$}}
\put(350,35){\makebox(0,0){$\g_1$}}
\put(400,35){\makebox(0,0){$\g_0$}}
\put(250,15){\makebox(0,0){$-2$}}
\put(300,15){\makebox(0,0){$-2$}}
\put(350,15){\makebox(0,0){$-2$}}
\put(400,15){\makebox(0,0){$-2$}}
\put(325,25){\makebox(0,0){$\cdots$}}
\end{picture}

\noindent 
Now, Theorem 5.4 follows from the above discussions and the next proposition.\\
5.14.\ {\bf Proposition.}\ {\em Assume that \x is numerical 
Gorenstein elliptic singularity with $p_g=m+1$ ($m>0$) and $\zn^2=-1$. 
Consider its minimal resolution graph and write $\g=B_1\cup \g_0$
(cf. 5.13). Then $|\co(-\zn)|$ has a unique basepoint $Q\in \g_0\setminus
B_1$, and $m_p\co_M=m_Q\co_M(-\zn)$. }\\
The {\em proof} is similar to \cite{Laufer77} (3.13), 
and it is left to the reader. \\

\noindent {\bf 6.\ The Hilbert--Samuel function and the embedding dimension.}\\

Assume that \x is an elliptic Gorenstein singularity with $p_g=m+1$.
Consider its minimal resolution and set $\cl:=\co_{\zk}(-\zn)$ and 
$d:=-\zn^2$. \\
6.1.\ {\bf Theorem.} \ {\em If $d\geq 3$ then $\oplus_{k\geq 0}
H^0(\zk,\cl^{\otimes k})$ is generated by elements of degree $k=1$.}\\
{\em Proof.}\ We will follow -- and modify -- the proof 
of the 
corresponding statement for the minimally elliptic singularities
(as it is presented in \cite{MR}, page 114;
i.e. ``Castelnuovo's free pencil trick'').

By (5.4) there exists $s_0\in H^0(\cl)$ with 
$m_{\g_i}(s_0)=m_{\g_i}(\zn)$ for all $\g_i$. Its divisor has form $\zn+D$,
where $D$ is an effective Cartier divisor with $|D|\cap \g$ finite, and
$\g_iD=\g_i\cl$. Obviously $\co_{\zk}(D)\to \cl$ (given by $f\mapsto
fs_0$) is an isomorphism. We use another generic section $s\in H^0(\cl)$
to identify $\cl$ with $\co_{\zk}$ near $D$, hence we will write $\co_D$
instead of $\cl|_D=\co_D(D)$. By construction, $\co_D$ is an Artinian 
scheme of degree $d$. Consider the exact sequence
$ 0\to \co_{\zk}\to \cl\to \co_D\to 0$, where the second arrow is the 
multiplication by $s_0$, (and the third is the ``division by $s$ near $D$''). 
Its long cohomology exact  sequence is:
$$0\to H^0(\co_{\zk})\stackrel{s_0}{\longrightarrow}H^0(\cl)\stackrel{r_1}
{\longrightarrow}
\co_D\to H^1(\co_{\zk})\to H^1(\cl)\to 0.$$
Since $h^1(\co_{\zk})=p_g$ and $h^1(\cl)=p_g-1$ (cf. 2.21.b), the rank of
$r_1$ is $d-1$. 

In the following lemmas $D'\subset D$ is an arbitrary subscheme of length 
$d-1$. Let $\ci\subset \co_{\zk}$ be its ideal sheaf. If ${\cal I}_D\subset 
\co_{\zk}$ is the ideal sheaf of $D$, 
then ${\cal I}_D\subset \ci$ and $\ci/{\cal I}_D=\bc_Q$
for some point $Q\in \g$. \\
6.2.\ {\bf Lemma.}\ $h^1(\ci\cl)=p_g-1$.\\
{\em Proof.}\ The inclusion $\ci\co_{\zk}(D)\subset \co_{\zk}
(D)=\cl$ guarantee that $\ci\cl$ is torsion free. Hence the exact sequence
\ $0\to \co_{\zk}\stackrel{s_0}{\to}\ci\cl\to\bc_Q\to 0$ \ does not split.
Hence in the exact sequence 
$0\to Hom(\ci\cl,\co_{\zk})\to Hom (\co_{\zk},\co_{\zk})\stackrel{\alpha}{\to}
Ext^1({\bf C}_Q,\co_{\zk})={\bf C}$ the map  $\alpha$ is not trivial. Hence by 
Serre duality $h^1(\ci\cl)=h^1(\co_{\zk})-1$. \ \ \ $\Box$\\
6.3.\ {\bf Lemma.}\ {\em $H^0(\cl)\to \co_{D'}$ is surjective for any
subscheme $D'\subset D$ of length $d-1$.}\\
{\em Proof.}\  Use the cohomology exact sequence of $0\to \ci\cl\to \cl\to
\co_{D'}\to 0$ and $h^1(\ci\cl)=h^1(\cl)$ (by 6.2. and 2.21.b). \ \ \ $\Box$\\
Now, by similar arguments as in \cite{MR}, pages 114-115,  Theorem 6.1 and
Theorem D (of the introduction) follows. \ \ \ $\Box$\\
6.4.\ {\bf Remarks.}

(a)\ {\bf The cases $d=-\zn^2=1,2$.}\ 
Similarly  to the minimally elliptic singularities
(cf. e.g. \cite{MR}, page 116), one  can write a generator set of the
$\co_X$--algebra $R:=\oplus_{k\geq 0}H^0(\co(-k\zn))$ in the case $d\leq 2$ as
well. 
If $d=2$, then $R$ is generated by $m_p$ in degree 1 and an element $y
\in m_p\setminus m_p^2$ in degree two.
If $d=1$, then $R$ is generated by $m_p$ in degree 1, by an element 
$y\in m_p\setminus m_p^2$ in degree 2, and $z\in m_p\setminus (m_p^2,y)$
in degree 3.

(b)\ Moreover (cf. \cite{Laufer77}), 
if $\zn^2=-4$, then \x is a (tangential) complete 
intersection; if $\zn^2\leq -5$, then \x is not a complete intersection.

(c)\ Using similar arguments as in this section, one can verify that
for any normal surface singularity $emb\, dim(X,p)\geq mult(X,p)+
min_{D>0}\chi(D)$. Actually, for any elliptic Gorenstein singularity
$emb\, dim(X,p)=mult(X,p)$ provided that $mult(X,p)\geq 3$, and 
(obviously) $emb\, dim(X,p)=3$ if $mult(X,p)=2$. In particular, for 
elliptic hypersurface singularities $mult(X,p)\leq 3$.

(d)\ \y computed the multiplicity and the Hilbert--Samuel function
for many different situations. His basic assumption was the much stronger
$Z_E^2\leq -2$ for the first case and 
$Z_E^2\leq -3$ for the second case. (For details see 
his series of papers listed in the References.)

(e)\ If $\zn^2=-1$, then (5.13) shows that the dual resolution graph is
a ``Kodaira graph'' (cf. e.g. \cite{K1}, 2.7). Moreover, by (5.14), the 
maximal ideal cycle is exactly $\zn$.  Hence, by \cite{K1} (2.9.1),
\x is  a Kodaira singularity. \\

\end{document}